\documentclass[12pt]{article}

\usepackage{graphicx}        % standard LaTeX graphics tool
\usepackage{multicol}        % used for the two-column index
\usepackage{amssymb}
\usepackage{amsmath, fancyhdr}
\usepackage{subfigure}
\usepackage{array}
\usepackage{rotating}
\usepackage{multirow}
\usepackage{booktabs}
\usepackage{tikz}
\usepackage{chemfig} 
\usepackage{dcolumn}
\usepackage{bm}
\usepackage[sort]{natbib}
\bibpunct{(}{)}{;}{a}{,}{,}
\usepackage[a4paper,left=1in,top=1in,width=6.3in,height=9in]{geometry}
\usepackage{authblk}

\newcommand{\nil}[1]{}

\newcolumntype{d}[1]{D{.}{\cdot}{#1}}
\newcolumntype{.}{D{.}{.}{-1}}
\newcolumntype{,}{D{.}{.}{4}}
\newcolumntype{;}{D{.}{.}{0}}
\newcolumntype{:}{D{.}{.}{2}}

\newcommand{\Kent}{\mbox{\text{Kent}}}
\renewcommand{\/}[1]{\bm #1}
\renewcommand{\geq}{\geqslant}
\renewcommand{\leq}{\leqslant}
\newcommand{\ignore}[1]{}
\newcommand{\kappahat}{\widehat{\kappa}}
\newcommand{\betahat}{\widehat{\beta}}
\newcommand{\gammahat}{\widehat{\/\gamma}}
\newcommand{\pihat}{\widehat{\pi}}
\newcommand{\acos}{\text{acos}}

\begin{document}

\title{Mixture models for spherical data with applications to
  protein bioinformatics}
\author[$^{1,2}$]{Kanti V.~Mardia}
\author[$^1$]{Stuart Barber}
\author[$^1$]{Philippa M.~Burdett}
\author[$^1$]{John T.~Kent}
\author[$^3$]{Thomas Hamelryck}
\affil[$$]{K.V.Mardia@leeds.ac.uk, mardia@stats.ox.ac.uk, S.Barber@leeds.ac.uk, J.T.Kent@leeds.ac.uk, thamelry@binf.ku.dk \vspace*{5mm}}
\affil[$^1$]{School of Mathematics, University of Leeds, LS2 9JT, UK}
\affil[$^2$]{Department   of Statistics, University of Oxford, OX1 3LB, UK}
\affil[$^3$]{Department of   Biology/Department of Computer Science,
  University of Copenhagen,   Denmark}

%%\authorrunning{Mardia, Barber, Burdett, Kent \& Hamelryck}% \emph{et al.}}
\date{\empty} \maketitle

\abstract{Finite mixture models are fitted to spherical data. Kent
  distributions are used for the components of the mixture because
  they allow considerable flexibility. Previous work on such mixtures
  has used an approximate maximum
  likelihood estimator for the parameters of a single component. However, the
  approximation causes problems when using the EM algorithm to
  estimate the parameters in a mixture model. Hence the exact
  maximum likelihood estimator is used here for the individual
  components.  \newline
  \hspace*{\parindent} This paper is motivated by a challenging prize
problem in structural bioinformatics of how proteins fold.  It is
known that hydrogen bonds play a key role in the folding of a
protein. We explore this hydrogen bond geometry using a data set
describing bonds between two amino acids in proteins.  An appropriate
coordinate system to represent the hydrogen bond geometry is proposed,
with each bond represented as a point on a sphere.  We fit mixtures of
Kent distributions to different subsets of the hydrogen bond data to
gain insight into how the secondary structure elements bond together,
since the distribution of hydrogen bonds depends on which secondary
structure elements are involved.}

%\abstract{Each chapter should be preceded by an abstract (no more than 200 words) that summarizes the content. The abstract will appear \textit{online} at \url{www.SpringerLink.com} and be available with unrestricted access. This allows unregistered users to read the abstract as a teaser for the complete chapter.\newline\indent
%Please use the 'starred' version of the \texttt{abstract} command for typesetting the text of the online abstracts (cf. source file of this chapter template \texttt{abstract}) and include them with the source files of your manuscript. Use the plain \texttt{abstract} command if the abstract is also to appear in the printed version of the book.}

\section{Introduction}
\label{sec:1}
Spherical data arise in many applications, such as palaeomagnetism
\citep{Hospers1955} optical crystallography
\citep[p.92]{Vistelius1966} and stellar directions \citep{Jupp1995};
more examples are given in \citet{Fisher1987} and \citet{Mardia2000}.
\cite{Pewsey2020} have given a recent review with many other
references.  In simple cases, these can be described as a sample from
an appropriate spherical distribution such as the Kent
distribution~\citep{Kent1982}.  However, many data sets are not so
straightforward and are better described as a mixture of components,
such as the rock fracture data considered by \citet{peel2001fitting}.
Often, we must estimate which observations come from each component
and a natural way to fit such a model is to use the EM algorithm.

\citet{Kent1982} developed a bivariate normal approximation for the
density of the Kent distribution.  This corresponding approximate
likelihood was used to construct a simple approximate maximum
likelihood estimate (MLE) of the parameters in place of the more
computationally intensive exact maximum likelihood estimator.
\citet{peel2001fitting} used the approximate likelihood in an EM
algorithm to fit a mixture of Kent distributions.  Unfortunately, with
this approximation the resulting EM algorithm is not necessarily
monotone (the EM algorithm must be monotone when the exact likelihood
is used).  The main problem is the use of the bivariate
normal approximation for the normalizing constant instead of the exact
value given below in (\ref{eq:norm-const}).  In this paper the exact
MLEs have been used to avoid the problem.

Our motivating example is the distribution of hydrogen bonds between
secondary structures in proteins.  Our data consist of four distances,
three angles, and two categorical variables. Driven by these
measurements and previous work by \cite{Paulsen2009} a new
`Euclidean-Latitude' representation of the hydrogen bond is
introduced, which reduces the original 12 variables to four
directional variables, incorporating three spherical coordinates and a
colatitude.  We focus on two of the angles which capture
the most important biological information about the hydrogen bond.

The joint distribution of these angles forms a `shell' on the surface
of the sphere, with the depth of the shell determined by the length of
the hydrogen bond. This distribution is dependent on the secondary
structures involved and the separation distance between the
interacting pair.  \cite{paulsenLASR} used bond lengths and PCA to
develop a probabilistic model for the bonds.  We propose a flexible
model for the distribution of hydrogen bonds by conditioning on the
hydrogen bond length and using the EM algorithm to fit a mixture
distribution of an unknown number of Kent distribution components
using the exact  MLEs of the parameters.  

\subsection{Proteins and hydrogen bonding}

For over forty years it has been known that a globular protein
spontaneously folds into a compact structure determined by its amino
acid sequence, but how the structure depends on the sequence is still
unknown.  A valuable tool in structure prediction is the ability to
simulate plausible candidate structures; examples are the methods
described by~\cite{HKK2006} and~\cite{ZPX2010}.  It is known that
hydrogen bonds play a key role in folding a protein into its 3-D
structure \citep{Baker1984}.  Each hydrogen bond can be defined by
four atomic coordinates in 3-D space in a representation equivalent to
that of \cite{kortemme2003orientation}.  We use an alternative
parameterisation based on directional statistics and focus on a
probabilistic approach to modelling hydrogen bonds.  Such a model has
potential applications in protein structure prediction, simulation,
validation and design. For example, the model could be used to
evaluate the potential strength of a hydrogen bond. By evaluating the
probability of the hydrogen bond network in a protein, potential
errors or low quality protein structures could be identified.

Proteins are polypeptide chains of amino acids joined end to end by
peptide bonds; see for example \cite{Branden} and \citet{Mardia2013}.
The sequence of these amino acids is known as the primary structure of
the protein.  There are twenty amino acids which all have a common
structure consisting of a central carbon atom, denoted $C_\alpha$, a
hydrogen atom, an amino group and a carboxyl group. A side chain,
unique to each of the 20~amino acid types, projects from the central
carbon atom.  As the polypeptide chain folds it forms secondary
structures which are classified into three types: helices, sheets and
loops.

\begin{figure}[!hbt]
  \begin{center}
    \includegraphics[width=0.5\textwidth]{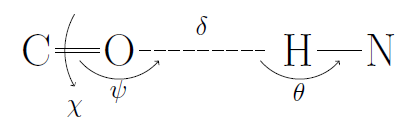}
    \caption{The three angles $\chi$, $\theta$, $\psi$ and 
      bond length $\delta$ associated with hydrogen bonds.  Note that
      the bonds C-O and H-N are assumed to be of constant length.}
    \label{fig:angles} 
  \end{center}
\end{figure}

The geometry of a hydrogen bond is shown in Figure \ref{fig:angles}.
A naive representation would be to use 12 variables (Euclidean
coordinates of the four atoms) to represent one bond.  However, the
spatial orientation and location of the bond are not considered
important, so we effectively lose two degrees of freedom for the fixed
bond lengths C-O and H-N, three for rotation and a further three for
translation.  Thus we retain full information about the geometry of
the hydrogen bond by using only four variables: the C-O-H angle
$\psi$, the dihedral angle N-C-O-H $\chi$ , the O-H distance $\delta$
and the N-H-O angle $\theta$ .  We focus on a subset of the angles,
$\chi$ and $\psi$.  These angles define a vector
$\/x = [x_1, x_2, x_3]^T$ on the surface of the unit sphere using
\begin{equation} \label{eq:x-psi-chi}
x_1 = \sin \psi \sin \chi, \quad x_2 = \sin \psi \cos
\chi,  \quad x_3 = \cos \psi, \quad \lvert \/x \rvert = 1,
\end{equation}
so that $\psi$ and $\chi$ represent the colatitude and longitude of
$\/x$.  Section~\ref{sec:eda} gives an analysis of the data using the angles
$\chi$ and $\psi$, together with the distance $\delta$, and a
separation distance $\Delta_L$ defined there.

\section{Fitting a single Kent distribution}

\subsection{The Kent distribution}

The Kent (or Fisher-Bingham 5 Parameter, $FB_5$) distribution
\citep{Kent1982} is the key distribution in the analysis of
directional data on $S^2$, the unit sphere in $\mathbb R^3$.  It is
the spherical analogue of the general multivariate normal distribution
allowing for distributions of any elliptical shape, size or
orientation on the surface of the sphere.  We write
$\Kent(\kappa, \beta, \/\Gamma)$ to denote the Kent distribution with
density
\begin{equation} \label{eq:kent-pdf}
  f(\/x) = \frac{1}{c(\kappa,\beta)}
  \exp\left\{ \kappa {\/\gamma_{(3)}^T\/x} +
  \beta [ (\/\gamma_{(1)}^T\/x)^2 - (\/\gamma_{(2)}^T \/x)^2]\right\}, 
  \quad |\/x| =1,
\end{equation}
where $\/\Gamma = [ \/\gamma_{(1)} ~ \/\gamma_{(2)} ~ \/\gamma_{(3)}]$
is a $3 \times 3$ orientation matrix and $c(\kappa, \beta)$ is a
normalizing constant.  Here, $\kappa$ determines the concentration of
the data, $\beta$ describes the ovalness of the contours of constant
probability, $\/\gamma_{(3)}$ is the population mean vector or pole of
the distribution, and $\/\gamma_{(1)},\/\gamma_{(2)}$ denote the major
and minor axes of the ellipse-like contours of constant probability.

An alternative angular parameterisation is obtained by noting that
$\/ \Gamma$ can be written as a product of three simple rotations,
$$
\Gamma = R_\phi R_\eta R_\omega,
$$
where
\begin{equation} \label{eq:rot-mats}
  \/R_\phi = \begin{bmatrix}
    \cos \phi & -\sin \phi & 0 \\
    \sin \phi & \cos \phi & 0 \\
    0 & 0 & 1
  \end{bmatrix}, \quad
 \/R_\eta = \begin{bmatrix}
    \cos \eta & 0 & \sin \eta \\
    0 & 1 & 0 \\
    -\sin \eta & 0 & \cos \eta \\
  \end{bmatrix}, \quad
   \/R_\omega = \begin{bmatrix}
    \cos \omega & -\sin \omega & 0 \\
    \sin \omega & \cos \omega & 0 \\
    0 & 0 & 1\\
  \end{bmatrix}.
\end{equation}
In particular, the final column of $\Gamma$ has elements
\begin{equation} \label{eq:x-theta-phi}
  \gamma_{13} = \sin \eta \sin \phi, \quad
  \gamma_{23} = \sin \eta \cos \phi, \quad
  \gamma_{33} =  \cos \eta,
\end{equation}
so that $\eta$ and $\phi$ represent the colatitude and latitude of the
mean direction.

Evaluation of the normalizing constant $c(\kappa, \beta)$ is not
trivial.  From \cite{Kent1982},
\begin{equation}
  \label{eq:norm-const}
  c(\kappa,\beta) = 2 \pi \sum_{j=0}^{\infty} \frac{\Gamma(j+ \frac{1}{2})}
  {\Gamma(j+ 1)}\beta^{2j}(\tfrac{1}{2}\kappa)^{-2j - \frac{1}{2}}
  I_{2j+ \frac{1}{2}}(\kappa);
\end{equation}
to calculate the normalizing constant we need to evaluate the modified
Bessel function $I_{2j+ \frac{1}{2}}(\kappa)$.  This can be calculated
using a backwards recursion which is more stable than forward
recursion.  From \cite{Amos},
\begin{align*}
I_{\nu - 1}(x) =& \frac{2\nu}{x}I_\nu(x) +I_{\nu+1}(x),\\
r_{\nu-1}(x)=&\frac{1}{2\nu/x+r_\nu(x)}, \quad r_\nu(x)=\frac{I_{\nu +1}(x)}{I_{\nu }(x)},\\
I_\nu(x)=& I_\alpha(x)\prod_{k=1}^{\lfloor\nu\rfloor}r_{\nu-k}(x),\quad \alpha = \nu-\lfloor\nu\rfloor, 
\end{align*}
where $\lfloor\nu\rfloor$ is the integer part of $\nu$.  Therefore,
\begin{align*}
I_{2j+ \frac{1}{2}}(\kappa)=&I_{\frac{1}{2}}(\kappa)\prod_{k=1}^{2j}r_{2j+\frac{1}{2}-k}(x),\\
\textrm{ where } I_{\frac{1}{2}}(\kappa)=&\sqrt{\frac{2}{\pi \kappa}}\sinh{\kappa} .
\end{align*}
The truncation point of the series, $m$, can be altered for the
desired level of accuracy.  In practice, we have found $m= 20$ to be
sufficient.

In passing, note that various approximations to the normalizing
constant have been developed over the years that improve over the
bivariate normal approximation in (\ref{eq:bvn}), notably the
saddlepoint approximation.  See, e.g., \citet{Kume-Wood2005} and
\cite{Scealy2014}.  However, the exact series used here
(\ref{eq:norm-const}) has the advantage that it can be computed to
arbitrary accuracy.

%%%%%%%%%%%%%%%%%%%%%%%%%%%%%%%%%%%%%%%%%%%%%

\subsection{High concentration bivariate normal approximation}

Under high concentration, a bivariate normal approximation to the Kent
distribution was given by \cite{Kent1982}.  Let $\/x\sim
\Kent(\kappa,\beta,\/\Gamma)$, and rotate to $\/y = \/\Gamma^T \/x$ so that
$\/y \sim \Kent(\kappa,\beta, \/I)$, where $\/I$ denotes the identity matrix.
The limiting distribution as $\kappa \rightarrow \infty$ and $\frac{\beta}{\kappa}
\rightarrow d$, where $0 \leq d < \frac{1}{2}$, can be described using
the bivariate normal distribution for the standardized variables
$\/y^*_1 = \sqrt{\kappa} y_1$ and $\/y^*_2 = \sqrt{\kappa} y_2$.  In particular
\begin{equation}
  \label{eq:bvn}
  \/y^*_1 \sim N\bigl(0,(1-2d)^{-1}\bigr), \quad \/y^*_2  \sim
N\bigl(0,(1+2d)^{-1}\bigr),
\end{equation}independently. The condition number of the corresponding
covariance matrix of the bivariate normal distribution is defined by
\begin{equation}\label{eq:eccentricity}
  \text{CN} = \frac{\kappa + 2 \beta} {\kappa - 2 \beta},
\end{equation}
and measures the extent of any anisotropy, where $d=0$
(i.e. $\beta=0$) for the isotropic case.

%%%%%%%%%%%%%%%%%%%%%%%%%%%%%%%%%%%%%%%%%%%%%

\subsection{Maximum likelihood estimators for the Kent distribution} \label{sec:Kent-mle}

Consider $n$ independent observations $\/x_1, \ldots,
\/x_n$ from a Kent distribution $\Kent(\kappa,\beta,\/\Gamma)$.  The log-likelihood 
is given by 
\begin{align}
  \ell(\Psi; \/x)
  &= \sum_{j=1}^n \left\{\kappa \/\gamma^T_{(3)} \/x_j +
    \beta [ (\/\gamma^T_{(1)}\/x_j)^2 - (\/\gamma^T_{(2)}\/x_j )^2] \right\} -
    n \log c(\kappa,\beta) \} \notag \\
  &= n  \kappa \/\gamma^T_{(3)} \bar{\/x} +
    n\beta \left(  \/\gamma^T_{(1)}\/S \/\gamma_{(1)} 
    - \/\gamma^T_{(2)}\/S\/\gamma_{(2)}\right)
    -  n \log c(\kappa,\beta) \label{eq:kent-log-lik}.
\end{align}
Here $\/\Psi=(\kappa,\beta, \eta, \phi, \omega)$ contains the parameters.
The log-likelihood depends on the data through the sample mean vector
and the sample second moment matrix
$$
\bar{\/x} = \frac1n \sum \/x_j, \quad
\/S = \frac1n \sum \/x_j \/x_j^T.
$$

The approximate MLE  of \citet{Kent1982} can be computed as follows.
\begin{itemize}
\item[(a)] The estimated orientation matrix $\hat{\/\Gamma}$ is
  determined so that $\hat{\/\Gamma}^T \bar{\/x} \propto [0\ 0\ 1]^T$
  and $(\hat{\/\Gamma}^T \/S \hat{\/\Gamma})_{12} = 0$, where the
  subscript indicates element $(1,2)$ of a $3 \times 3$ matrix.  This
  estimate depends just on the first two sample moments and does not
  require any estimates of $\kappa$ and $\beta$.

\item[(b)] The normalizing constant is approximated using the bivariate normal
  approximation in (\ref{eq:bvn}).  With this substitution the estimation of
  $\kappa$ and $\beta$ is straightforward.
  \end{itemize}

  Mixture models are discussed below in Section \ref{sec:kent-mixtures}.
  As noted earlier, the use of the approximate MLE in mixture models
  can cause problems for the EM algorithm, because the approximate
  log-likelihood can sometimes fail to increase at each iteration.
  Hence for this paper we limit attention to the exact MLE for the 
  single-component Kent distribution.  The exact MLE needs to be found
  numerically, e.g. by using the function \textsf{optim} in
  \textsf{R}, and for this purpose an analytic form for the
  derivatives of $\ell$ is very helpful.  For the concentration
  parameters
  \begin{align*}
    \partial \ell / \partial \kappa
    & =  n\/\gamma_{(3)}^T \bar{\/x} 
      -n \partial \log c(\kappa,\beta) /
      \partial \kappa, \\
    \partial l / \partial \beta
    & =   n \left(  \/\gamma^T_{(1)}\/S\/\gamma_{(1)}
      - \/\gamma^T_{(2)}\/S\/\gamma_{(2)} \right)
      -n \partial \log c(\kappa,\beta) / \partial \beta,
  \end{align*}
  where the derivatives of $c(\kappa)$ are obtained by differentiating
  the series (\ref{eq:norm-const}) termwise.  For the orientation
  parameters, the derivatives of $\Gamma$ are obtained by
  differentiating each of its building  blocks.  For example,
  $$
  \partial \/\Gamma/\partial \phi =  R'_\phi R_\eta R_\omega, \quad
  \/R'_\phi = \begin{bmatrix}
    -\sin \phi & -\cos \phi & 0 \\
    \cos \phi & -\sin \phi & 0 \\
    0 & 0 & 1
  \end{bmatrix},
  $$
  and similarly for $\partial \/\Gamma/\partial \eta$ and
  $\partial \/\Gamma/\partial \omega$.

%%%%%%%%%%%%%%%%%%%%%%%%%%%%%%%%%%%%%%%%%%%%%

\section{Fitting Kent mixtures using the EM algorithm}
\label{sec:kent-mixtures}
\subsection{The mixture distribution}

Consider a mixture of $g$ Kent components with density
\begin{equation} \label{eq:mix-pdf}
f(\/x;\/\Psi)=\sum_{i=1}^g \pi_i f(\/x; \/\Psi_i) 
\end{equation}
Here 
$\pi = (\pi_1, \ldots \pi_{g}$, where $\pi>0, \ i=1, \ldots, g$ and
$\sum \pi_i = 1$ is a vector of \emph{mixing proportions}.  
Also,  $f(\/x; \/\Psi_i)$ is the Kent density for the $i$th component,
$i=1, \ldots, g$ with parameter vector $\Psi_i$.  In many cases it is
useful to let one of the components denote the uniform distribution on
the sphere with density $f_0(\/x) = 1/(4 \pi)$.  Of course for the
uniform component there are no parameters to estimate.  This is the
same model used by \cite{peel2001fitting}.

Let $\/Z=(z_1,\ldots,z_n)^T$ be an $n \times g$ matrix membership
matrix where $z_{ji} =1$ if observation $j$ belongs to the $ith$
component and 0 otherwise.  In general the membership matrix is
unobserved and represents a set of missing or latent variables.
However, if the membership matrix is available, then a complete-data
vector $\/w_j$ for observation $j$ can be defined by
$\/w_j = (\/x_j^T, \/z_j^T)^T$.  The complete-data log-likelihood
function for $\/\Psi$ and $\/\pi$ is given by
\begin{equation} \label{eqn:totallikelihood}
  \ell_c(\/\Psi, \/\pi) =
  \sum_{i=1}^{g} \sum_{j=1}^n z_{ji} \{\log \pi_i + \log f(x_j; \/\Psi_i)\}.
\end{equation}

When the membership matrix is not observed, then it is possible to compute the
probability that each observation belongs to each group.  
Given parameter values   $\/\Psi_i, \ i=1, \ldots, g$ and group probabilities
$\pi_i, \ i=1, \ldots, g$, the
membership probabilities  for observation $j$  are given by 
\begin{equation} 
\label{eq:tau}
\tau_{ji} = P(Z_{ji}=1 | \/x_j) = \frac{\pi_i f(\/x _j; \/\Psi_i)}
{\sum _{h=1}^g \pi_h f(\/x _j; \/\Psi_h)}. 
\quad i=1, \ldots, g.
\end{equation}

%%%%%%%%%%%%%%%%%%%%%%%%%%%%%%%%%%%%%%%%%%%%%

\subsection{EM algorithm}

The EM algorithm \citep{dempster1977maximum} can be used to
iteratively compute the exact maximum likelihood estimator for a
mixture of $g$ Kent components.  It consists of alternating
expectation (E) steps and maximisation (M) steps.  It requires a
knowledge of the number of groups $g$ and initial estimates for the
parameters $\/\Psi_i^{(0)}, \ \pi_i^{(0)}, \ i=1, \ldots, g$.  Let
$\nu \geq 0$ index the current iteration.
\begin{itemize}
\item[(a)] (E-step). Given current values for the parameters
  $\/\Psi_i^{(\nu)}, \ \pi_i^{(\nu)}, \ i=1, \ldots, g$, estimate the
  membership probabilities using (\ref{eq:tau}) to get estimates
  $\tau_{ji}^{(\nu)}$, say.

\item[(b)] (M-step). Given the the estimated membership probabilities
  $\tau_{ji}^{(\nu)}$, compute weighted sample moments
  \begin{equation}
    \bar{\/x}_i^{(\nu)} = \frac{\sum_{j=1}^n \tau^{(\nu)}_{ji} \/x_j}
    {\sum_{j=1}^n \tau^{(\nu)}_{ji}}, \quad
    \/S_i^{(\nu)} = \frac{\sum_{j=1}^n \tau^{(\nu)}_{ji} \/x_j\/x_j^T}
    {\sum_{j=1}^n \tau^{(\nu)}_{ji}}.
  \end{equation}
  Using these weighted moments in (\ref{eq:kent-log-lik}), update the
  maximum likelihood estimates for the parameters
  $\/\Psi_i^{(\nu+1)}, \ i=1, \ldots, g$.  In addition update the
  estimates of the mixture proportions to
  \begin{equation}
    \label{eq:pis}
    \pi^{(\nu+1)}_i = \sum_{j=1}^n \tau^{(\nu)}_{ji}/n, \quad i=1, \ldots, g.
  \end{equation}
\end{itemize}
Then cycle between these two steps until convergence.  For more
details on the EM algorithm, see e.g.  \cite{little2002statistical}
and \cite{mclachlan2008algorithm}.

To select the number $g$ of components in the model, we choose that value
which minimises the AIC,
$$
2 k^* - \hat{\ell}_g.
$$
Here $g$ is the number of components in the model, and $\hat{\ell}_g$ denotes the maximized log-likelihood for the $g$-component mixture.  If all the components
are Kent distributions, then the total number of parameters is
$k^* = 5g + g-1 = 6g-1$ (5 parameters for each Kent distribution and $g-1$
parameters for the mixing probabilities).  If one of the components is a uniform
distribution, then the total number of parameters is
$k^* = 4g + g-1 = 5g-1$.

%%%%%%%%%%%%%%%%%%%%%%%%%%%%%%%%%%%%%%%%%%%%%

\section{Simulation study}

To assess the performance of our EM algorithm for fitting the Kent
mixture model given by~(\ref{eq:mix-pdf}), we consider several
simulation experiments.  In section~\ref{sec:angular} we consider how
classification rate improves as components move further apart on the
surface of the sphere, assuming the number of components in known.  In
section~\ref{sec:select}, we use a range of conditions to examine the
accuracy of the method when the number of components must be estimated
by minimizing the AIC.

\subsection{Angular separation of components} \label{sec:angular}

Consider an observation $j$.  A set of classification probabilities
$\hat{\tau}_{ji}, \ i=1, \ldots, g$ can be \emph{hardened} by
assigning observation $j$ to the group $i$ with the largest
classification probability.  In any simulation experiment, the labels
of all the observations are known.  Hence for each true group $i$, the
number misclassified observations can be counted.  The results can be
summarized either as a misclassification rate for each group $i$ or as
an overall misclassification rate.

The misclassification rate gives a useful indication of how distinct the
groups are from one another. Clearly the angular separation of the
components will affect the misclassification rate, where the angular
separation between two unit vectors $\/l_1$ and $\/l_2$ is
measured by $\acos(\/l_1^T \/l_2)$.

To explore the misclassification rate, four Kent components each of
size $n=200$, and all with $\kappa=50$ and $\beta=10$ were generated.
One component was located at the north pole, one at the south pole and
two on the equator. The two equatorial components ($E_1$ and $E_2$)
initially had the same mean direction. The two polar components, and
$E_1$ were kept stationary while $E_2$ was rotated around the equator
so the angular separation ranged from $0$ to $\pi$.  A 4-component
Kent model (without a uniform component) was fitted to the data using
the known values for the mean direction, $\kappa$, $\beta$, and
$\pi_i$ as starting points.

\begin{figure}[!hbt]
\centering
 \includegraphics[width=0.5\textwidth]{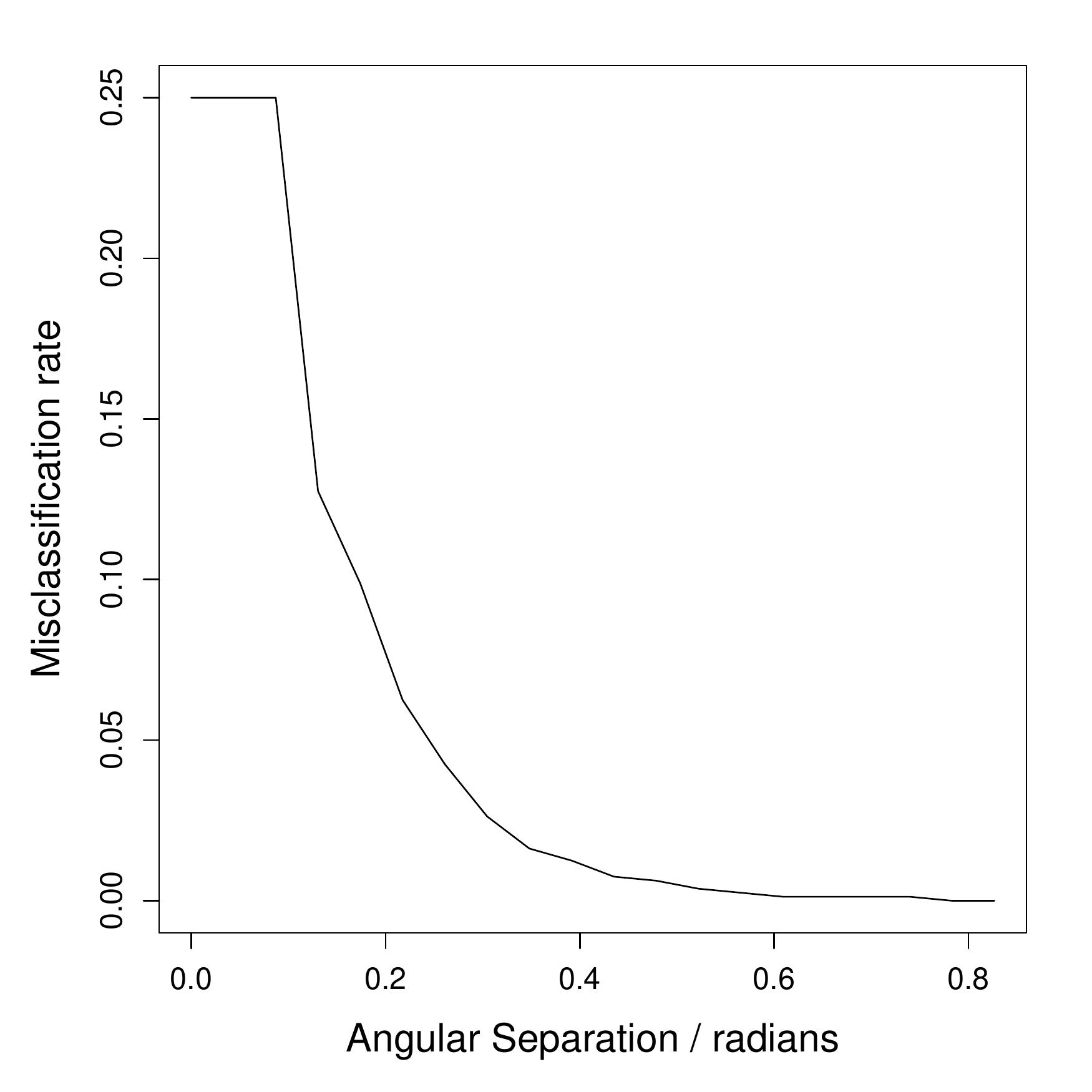}
 \caption[Misclassification rate]{Misclassification rate against
   angular separation $\acos(E_1^T E_2)$.}
\label{fig:misclass2}
\end{figure}

The misclassification rate as the angular separation between $E_1$ and
$E_2$ increased is shown in Figure \ref{fig:misclass2}.  As expected,
the misclassification decreases as the angular separation increases
and the overlap between $E_1$ and $E_2$ decreases.
 
%%%%%%%%%%%%%%%%%%%%%%%%%%%%%%%%%%%%%%%%%%%%%

\subsection{Model selection} \label{sec:select}

We shall illustrate our method with data simulated under three
conditions.  In Case 1, four samples of equal size, $n=200$, were
generated from the Kent distribution with different values for
$\kappa, \beta$ and the mean vector $\/\gamma_{(3)}$.  Case 2 looks at
five samples of differing sizes and values for the parameters.  In
Case 3 the Case 2 data set is re-used with 100 spherical uniform
random variates added.  In all cases a uniform component is included in the
mixture model fitted to the data.

%%%%%%%%%%%%%%%%%%%%%%%%%%%%%%%%%%%%%%%%%%%%%

\subsection{Case 1: Four equi-sized components}

In this case samples of size $n=200$ were simulated from $g=4$  Kent distributions
with parameters given in Table \ref{tab:case1}. Mixture models defined
by~(\ref{eq:mix-pdf}) with $1, 2, \ldots, 7$ Kent components, plus a
uniform component, were fitted to the data set and the AIC for each model was computed.

\begin{table}[t]
\caption{Case 1. True simulation parameters (top panel) and estimates
  of these parameters for the model with the lowest AIC (lower
  panel). }
\label{tab:case1}
\begin{center}
\begin{tabular}{cccccc}
\toprule
 & \multicolumn{5}{c}{Component} \\ \cline{2-6}
 & uniform& 1 & 2 & 3 & 4\\
\midrule
$\kappa$ & - & 5 & 10 & 20 & 20\\[1mm]
$\beta$ & - & 2 &4 & 9 & 2\\[1mm]
$\pi_i$ & 0 &0.25& 0.25& 0.25& 0.25\\[1mm]
$\/{\gamma_{(3)}}^T$  & - & $[0,0,1]$ & $[0,1,0]$ &$[0,-1,0]$ & $[0,0,-1]$ \\
\midrule
$\kappahat$ &- &6.94 &10.70 &18.22 &22.16\\[1mm]
$\betahat$& - &3.14& 5.00& 8.55& 2.82\\[1mm]
$\pihat_i$& 0.06 &0.21& 0.24 &0.24& 0.24\\[1mm]
$\/{\gammahat_{(3)}}^T$ & -	
	& $[0.0, 0.0, 1.0]$
	& $[0.0, 1.0, 0.0]$
	& $[0.0, -1.0, 0.0]$ 
	& $[0.0, 0.0, -1.0]$\\	
\bottomrule
\end{tabular}
\end{center} 
\end{table}

The fitted model  with the smallest AIC has the
correct choice of four Kent components (plus a uniform component) and
provides MLEs that are very close to the true
parameters.  In addition the misclassification rate is $2.5\%$, with
eleven observations allocated to the uniform component.  In this case,
the mean directions of the components are well-separated.

%%%%%%%%%%%%%%%%%%%%%%%%%%%%%%%

\subsection{Case 2: Five unequally sized components}

In this case $n=200$ samples were simulated from $g=5$ Kent
distributions with parameters given in Table \ref{tab:case2}.  Mixture
models with $1,\ldots,8$ Kent components, plus a uniform component,
were fitted to the data set, and the AIC for each model was computed.
The model with the smallest AIC has the correct choice of five Kent
components.  The parameter estimates for this model are shown in Table
\ref{tab:case2}

\begin{figure}[tbh]
\centering
 \includegraphics[ width=\textwidth]{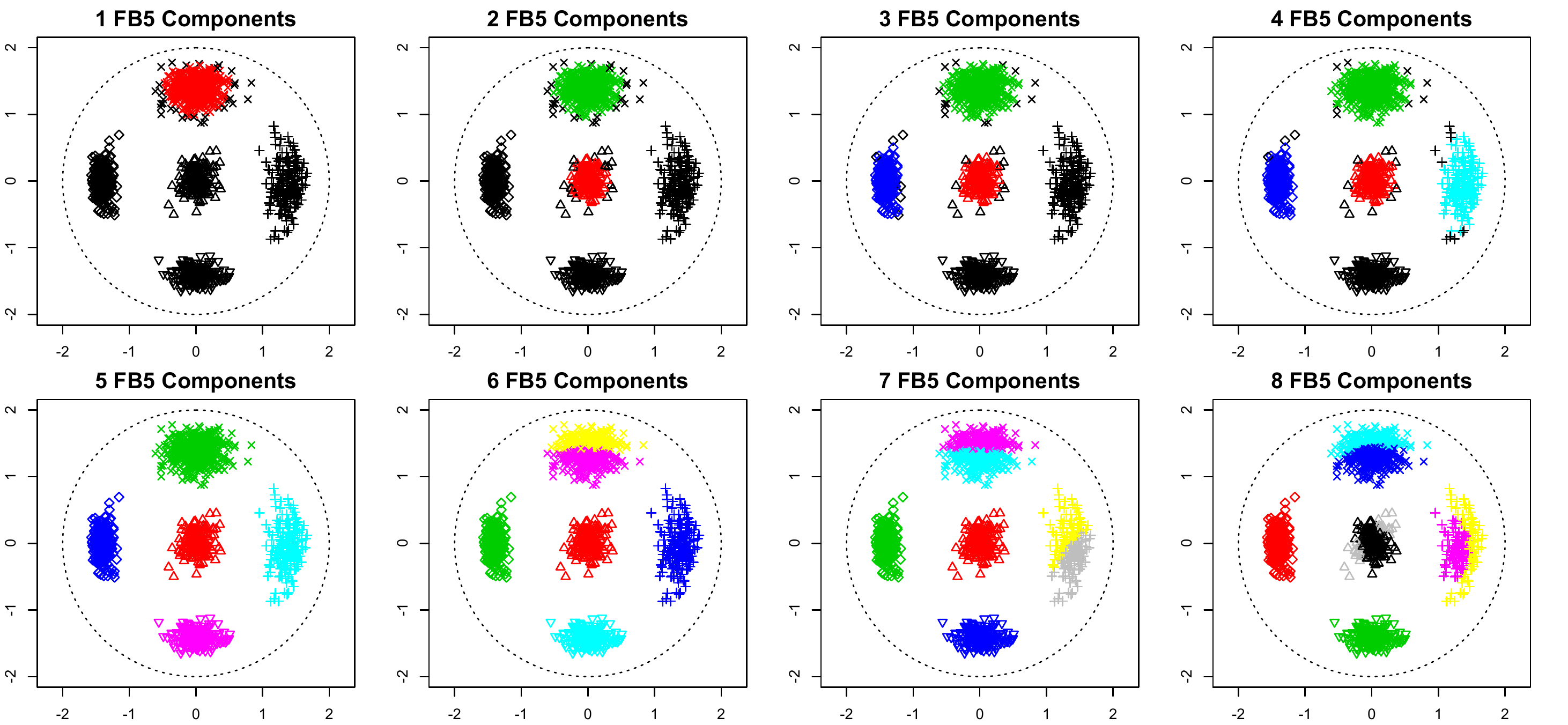}
 \caption{Schmidt net projections of the data coloured according to
   each model.  The plotting symbols indicate the known membership of
   the data.  The top left panel shows the data with only one Kent
   component fitted, plus a uniform component. In each subsequent
   panel an additional Kent component has been added to the model.
   The fitted model with five Kent components (the correct number) is
   shown as the first plot in the bottom row.  Points allocated to the
   uniform component are in black in each panel.}\label{fig:sn}
\end{figure}

Figure \ref{fig:sn} shows the Schmidt-net projections of the data
coloured by their allocations for each of the models fitted.  The
plotting symbols show the known allocations.  For Schmidt-net
projection, see, for example, \citet[p.161]{Mardia2000}; this is an equal
area projection of the unit sphere in $\mathbb R ^3$ onto a disc of
radius two where
$$(y_1,y_2)=2\sin (\psi / 2) (\cos \chi, \sin \chi), $$ where the
polar coordinates $(\psi, \chi)$ are defined in~(\ref{eq:x-psi-chi}).
From these plots it can also be seen that the model with five Kent
components gives the best fit with no error, i.e.\ a misclassification
rate of $0\%$.  In particular, although a uniform component has been
included in the fitted model, no points have been allocated to it.

With fewer than five components in the model, some of the true components are
entirely allocated to the uniform component,  and for those components
that are fitted correctly,  a few observations have also been allocated to be
uniform --- in these cases the fitted $\kappa$ values are overly large
making the fitted model more concentrated.  For the model with more
than five Kent components the clusters have been split into two or
more components.  Again we see that the overly complex model results in
some observations being erroneously assigned to the uniform component.

\begin{table}
\caption{Case 2: True simulation parameters (top panel) and estimates
  (rounded to 2 d.p. to highlight the effect on  $\/\gamma_{(3)}$) for
  the selected model with the lowest AIC (lower
  panel).} \label{tab:case2}
\begin{center}
\begin{tabular}{c,,,,,}
\toprule
&\multicolumn{5}{c}{Component}\\ \cline{2-6}
 & 1 & 2 & 3 & 4& 5\\
\midrule 
$\kappa$ &50& 25 & 30 & 70 & 50\\[1mm]
$\beta$ & 2& 2 &5& 10	&2\\[1mm]
$\pi_i$ &0.15& 0.15 & 0.37& 0.19& 0.15\\[1mm]
$\/{\gamma_{(3)}}$ &
\multicolumn{1}{c}{$\left[\begin{array}{c}0\\0\\1 \end{array}\right]$}& 
\multicolumn{1}{c}{$\left[\begin{array}{c}0\\1\\0 \end{array}\right]$}& 
\multicolumn{1}{c}{$\left[\begin{array}{c}1\\0\\0 \end{array}\right]$}& 
\multicolumn{1}{c}{$\left[\begin{array}{c}~~0\\-1\\~~0 \end{array}\right]$}& 
\multicolumn{1}{c}{$\left[\begin{array}{c}-1\\~~0\\~~0 \end{array}\right]$}
  \\[2mm]
\midrule
$\kappahat$ & 50.07& 24.97& 29.54& 68.22 &52.31\\[1mm]
$\betahat$ & 5.76&  2.71&  5.93& 11.08 & 4.55\\[1mm]
$\pihat_i$ & 0.15 &0.15 &0.37& 0.19& 0.15\\[1mm]
$\gammahat_{(3)}$ & 
\multicolumn{1}{c}{$\left[\begin{array}{r} -0.02\\ 0.01\\ 1.00\end{array} \right]$}& 
\multicolumn{1}{c}{$\left[\begin{array}{r} 0.00\\ 1.00\\ 0.04\end{array} \right]$}& 
\multicolumn{1}{c}{$\left[\begin{array}{r} 1.00\\ -0.01\\ 0.00\end{array} \right]$}&
\multicolumn{1}{c}{$\left[\begin{array}{r} 0.00\\ -1.00\\ 0.00\end{array} \right]$}& 
\multicolumn{1}{c}{$\left[\begin{array}{r} -1.00\\ -0.02\\ 0.00\end{array} \right]$}\\
\bottomrule
\end{tabular}
\end{center}
\end{table}

%%%%%%%%%%%%%%%%%%%%%%%%%%%%%%%%%%%%%%%%%%%%%

\subsection{Case 3: Case 2 with a uniform component added}

To the previous data set (Case 2) of five Kent components, we added
100 random uniform variates. The algorithm was then applied with the
same starting points.  In this case, the AIC picks out the correct
model with five Kent components (plus a uniform component), with a
misclassification rate of $2.6\%$.  Table~\ref{tab:misclass} shows
that over $75\%$ of the misclassified points were from the uniform
component and of these all can be shown to be close to a Kent
component.  Hence we conclude that our fitting procedure for the Kent
mixtures performs well for components with little overlap.

\begin{table}
  \caption{Case 3. Estimates of the parameters (rounded to 2 d.p.\ to
    show the effect on $\/\gamma_{(3)}$) when a uniform component is
    included.  The true parameters are shown in Table
    \ref{tab:case2}.\label{tab:punif}}
\begin{center}
\begin{tabular}{c,,,,,,}
\toprule
&\multicolumn{6}{c}{Component}\\ \cline{2-7}
 & \textrm{uniform} & 1 & 2 & 3 & 4& 5\\
\midrule  
$\kappahat$&\multicolumn{1}{c}{--}& 56.01& 24.67& 29.79 &67.48 &51.23\\[1mm]
$\betahat$&\multicolumn{1}{c}{--}&  5.53&  2.62&  6.21& 10.17 & 4.18\\[1mm]
$\pi_i$ &0.07& 0.14& 0.14& 0.34& 0.17& 0.14\\[1mm]
$\pihat_i$& 0.08& 0.13 &0.14& 0.35& 0.17& 0.14\\[1mm]
$\gammahat_{(3)}$ & 
\multicolumn{1}{c}{-} &
\multicolumn{1}{c}{$\left[\begin{array}{r}-0.02\\ 0.01\\ 1.00\end{array}\right]$}& 
\multicolumn{1}{c}{$\left[\begin{array}{r}-0.01\\ 1.00\\ 0.03\end{array}\right]$}& 
\multicolumn{1}{c}{$\left[\begin{array}{r} 1.00\\-0.01\\ 0.00\end{array}\right]$}&
\multicolumn{1}{c}{$\left[\begin{array}{r} 0.00\\-1.00\\ 0.00\end{array}\right]$}& 
\multicolumn{1}{c}{$\left[\begin{array}{r}-1.00\\-0.02\\-0.01\end{array}\right]$}\\
\bottomrule
\end{tabular}
\end{center}
\end{table}

\begin{table}[!htb]
  \caption{Case 3: The classifications of the data points from the
    fitted model compared to the known
    classifications.} \label{tab:misclass} 
\begin{center}
\begin{tabular}{c@{~~~~}cccccc}
\toprule
 & \multicolumn{6}{c}{True component}\\ \cline{2-7}
Classification & Uniform &   Kent 1 &  Kent  2&  Kent   3 & Kent  4 &  Kent  5 \\ \midrule
  Uniform	& 71 	&   4	&   1&   2&   1&   1\\
  Kent 1 	&  1 	&196	& 0   &0   &0   &0\\
  Kent 2 	&  5	& 0 	& 199&  0 &  0 &  0\\ 
  Kent 3 	& 12 	&  0 	&  0& 498 &  0  & 0\\
  Kent 4 	&  6 	& 0  	& 0 &  0 &249   &0\\
  Kent 5	&  5 	& 0 	&  0 &  0  & 0 &199\\ \bottomrule
\end{tabular}
\end{center}
\end{table}

%%%%%%%%%%%%%%%%%%%%%%%%%%%%%%%%%%%%%%%%%%%%%

\section{Hydrogen bond data} \label{sec:hbond}

\subsection{Exploratory data analysis} \label{sec:eda}

We consider a data set of $11653$ potential hydrogen bonds.  It is
important to note that the hydrogen bonds in this data set are not
observed, but that the presence of a hydrogen bond is inferred between
two non-adjacent residues within the protein chain where the
oxygen-nitrogen bond has length less than $3.5 \mathring{A}$.  If
several connections satisfied these criteria, then the one with the
smallest oxygen-hydrogen distance, $\delta$, was selected, see
Figure~\ref{fig:angles}.  We therefore assume that some of the alleged
hydrogen bonds in the data are spurious, and could bias parameter
estimates.

The data on each potential hydrogen bond consists of angles $\chi$ and
$\psi$, as defined in~(\ref{eq:x-psi-chi}), separation distance
($\Delta_L$), and secondary structure combination.  Thus, our data are
in the coordinate system $(x_1, x_2, x_3)$ defined by $(\chi, \psi)$
in~(\ref{eq:x-psi-chi}).  The variable $\Delta_L$ is given by the
number of amino acid residues separating the interacting pair where
$\Delta_L \in \{2,3,4,5+\}$ as hydrogen bonds cannot form between
adjacent amino acids.

The data also contain information on the secondary structure of the
residues in the chain on the O-side and the N-side of the hydrogen
bond (see \cite{Branden}).  There are three types of secondary
structure ($\alpha$-helix, loop and $\beta$-sheet), and hydrogen bonds
cannot form between $\beta$-sheets and $\alpha$-helices, so there are
seven possible secondary structure combinations.  In this paper, we
restrict ourselves to modelling the $n = 3627$ potential hydrogen
bonds formed between $\alpha$-helices.

\begin{figure}[t]
  \begin{center}
    \begin{tabular}{ccc}
      \includegraphics[width=0.3\textwidth]{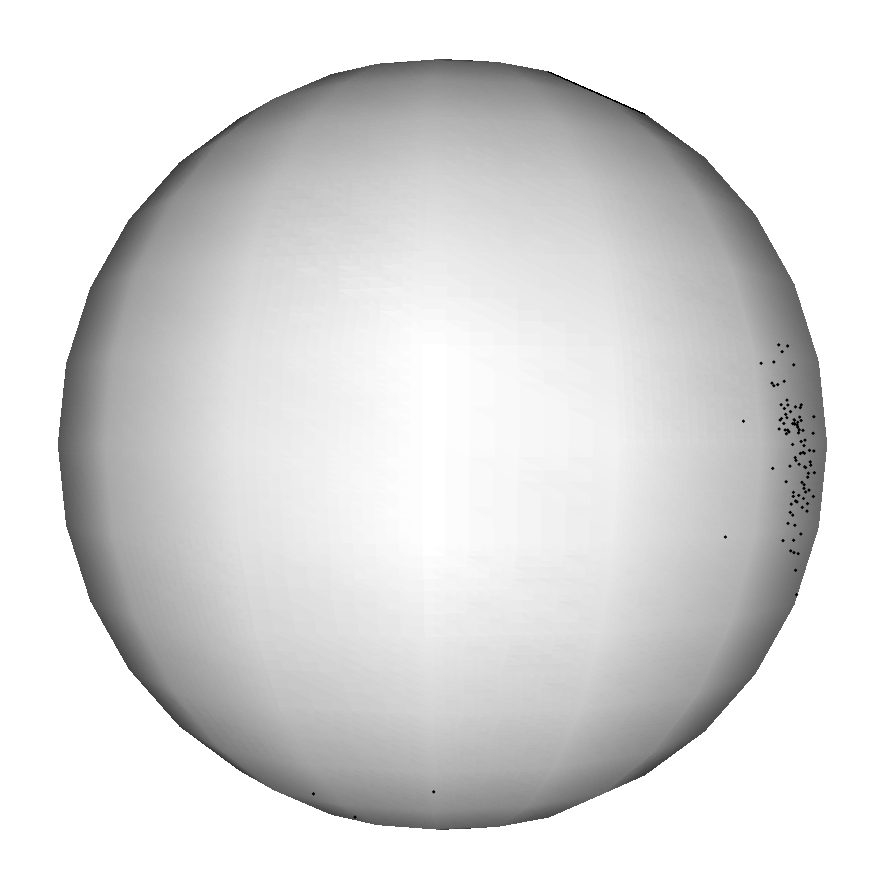} & 
      \includegraphics[width=0.3\textwidth]{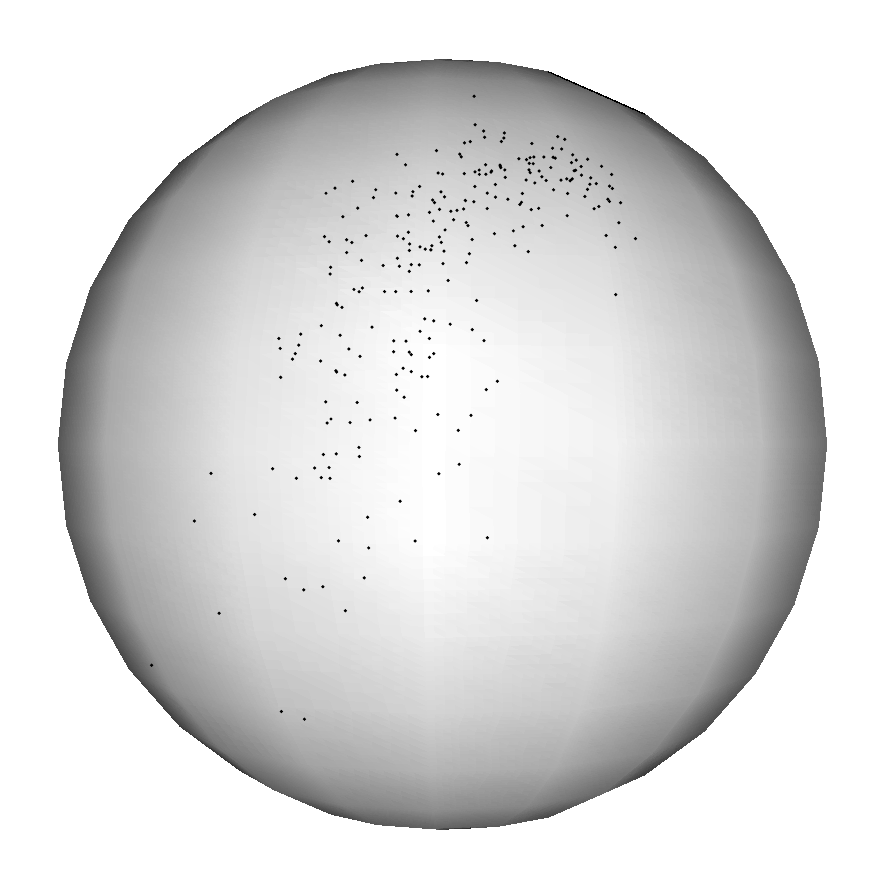} & 
      \includegraphics[width=0.3\textwidth]{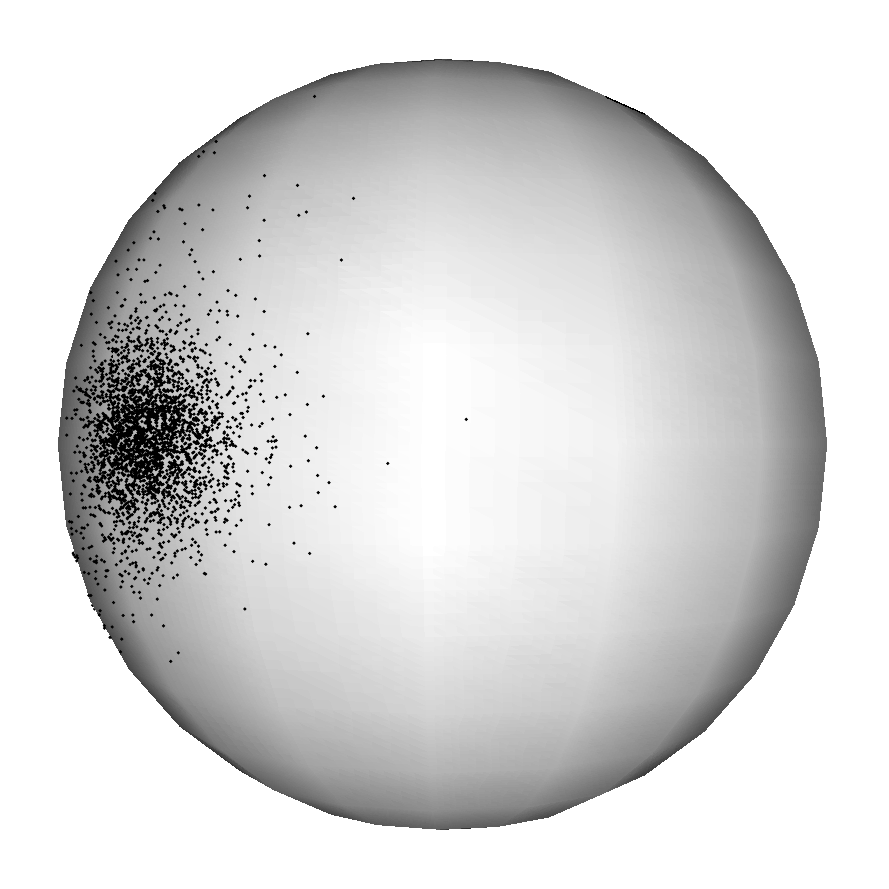} \\
      (a) & (b) & (c) \\
    \end{tabular}
    \caption{Plots showing the helix-helix data separated by
      $\Delta_L$, all shown from the same point of view.  Panels
      (a)--(c) show data for $\Delta_L=2, 3$, and 4 respectively.}
    \label{fig:HHspheres} 
  \end{center}
\end{figure}

Plots of the data can be made by projecting the data on to the sphere.
Using the separation distance $\Delta_L$ to break down the data
further showed that this distance affects the distribution of the
hydrogen bonds.  An example is shown in Figure \ref{fig:HHspheres} of
the helix-helix data, with most of the data having $\Delta_L=4$.  From
a biological point of view, it is to be expected that nearly all
hydrogen bonds between helical amino acids have $\Delta_L = 4$, since
this type of hydrogen bond provides the stability of the helix. The
remaining hydrogen bonds can be regarded as outliers or examples of
unusual hydrogen bond geometries.

%%%%%%%%%%%%%%%%%%%%%%%%%%%%%%%%%%%%%%%%%%%%%

\subsection{ Modelling the helix-helix data}

\begin{figure}[!hbt]
  \begin{center}
    \begin{tabular}{ccc}
      \includegraphics[ width=.49\textwidth]{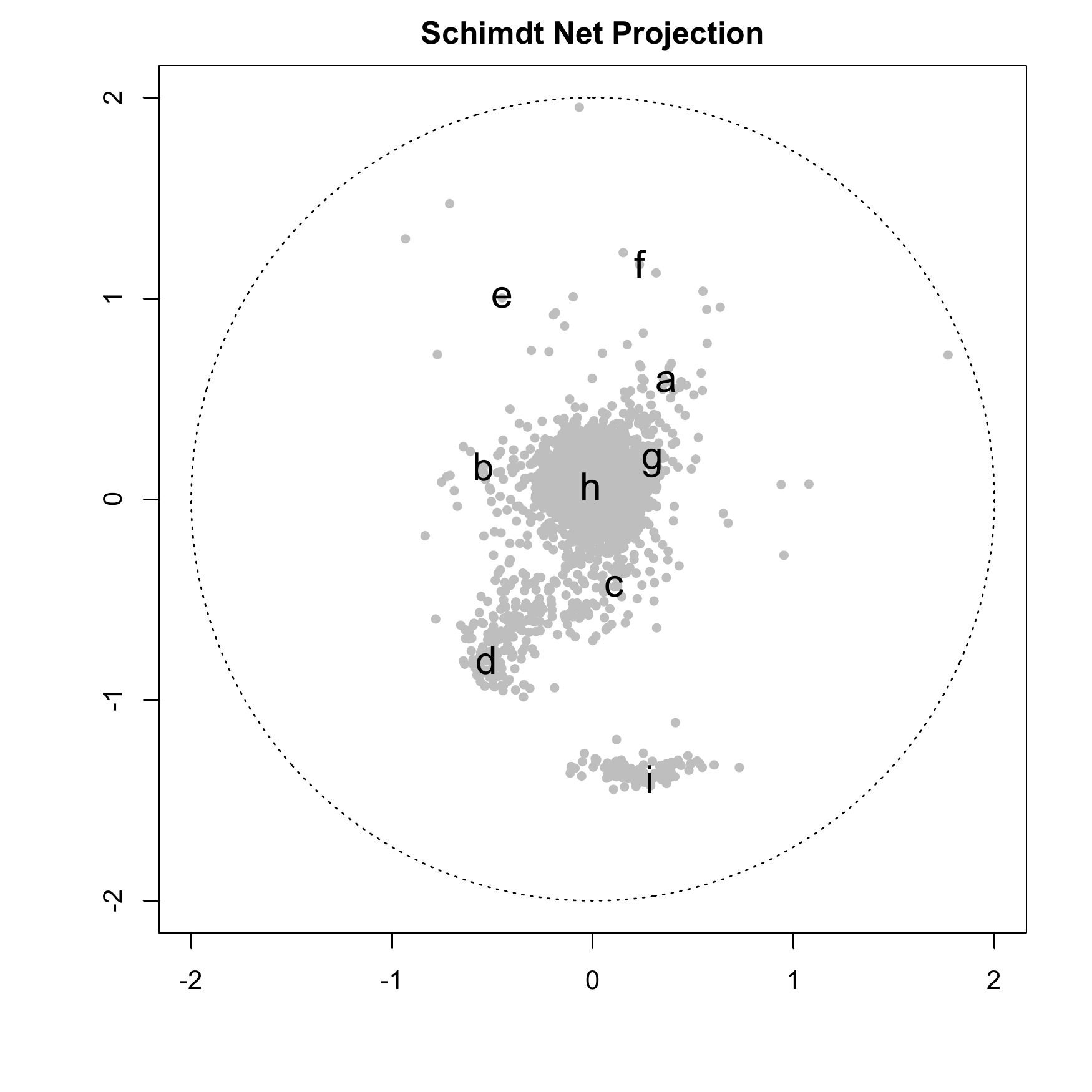} &
      \includegraphics[ width=.49\textwidth]{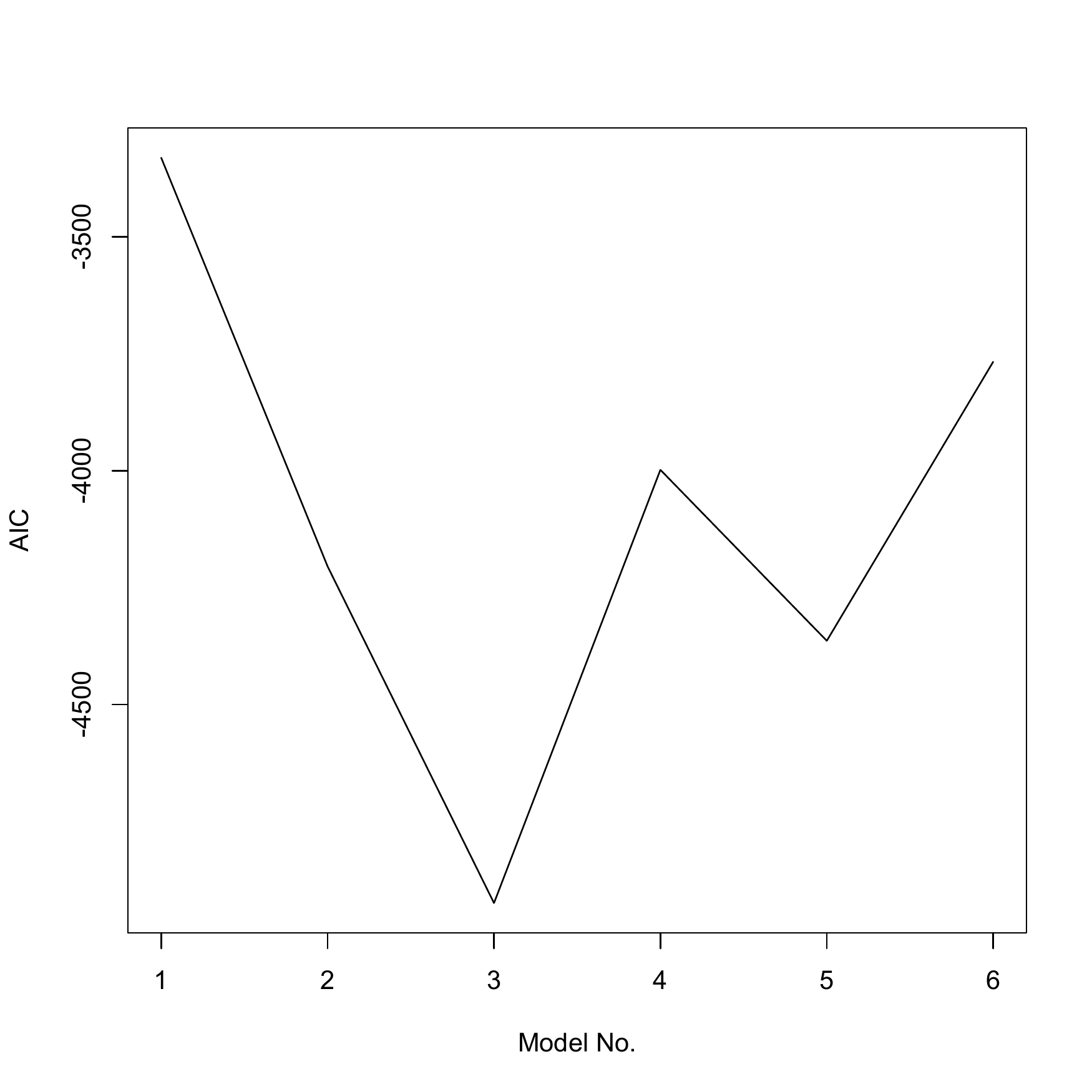} \\ 
      (a) & (b) \\
    \end{tabular}
  \end{center}
  \caption{Plots showing (a) Schmidt-net projection of the helix-helix
    data with eight different starting points $a,\ldots,i$ and (b) the
    AIC for each model fitted.  The AIC is minimised by the model with
    three Kent components.} \label{fig:hh}
\end{figure}

The helix-helix data shows one large cluster containing most of the
data, and several diffuse areas.  Choosing starting points for the
algorithm by conditioning on $\Delta_L$ works for the larger clusters
but cannot separate out the side clusters effectively.  Also this
method does not allow for choosing further starting points.  Therefore
the starting points for the mean directions were selected by plotting
a Schmidt-net projection of the complete helix-helix data (Figure
\ref{fig:hh}(a)) and identifying individual points by eye as possible
choices. Initial values for $\kappa, \beta$ and the mixing
proportions, $\pi_i$, were also chosen by visual inspection.

\begin{figure}[!hbt]
\centering
 \includegraphics[ width=0.7 \textwidth]{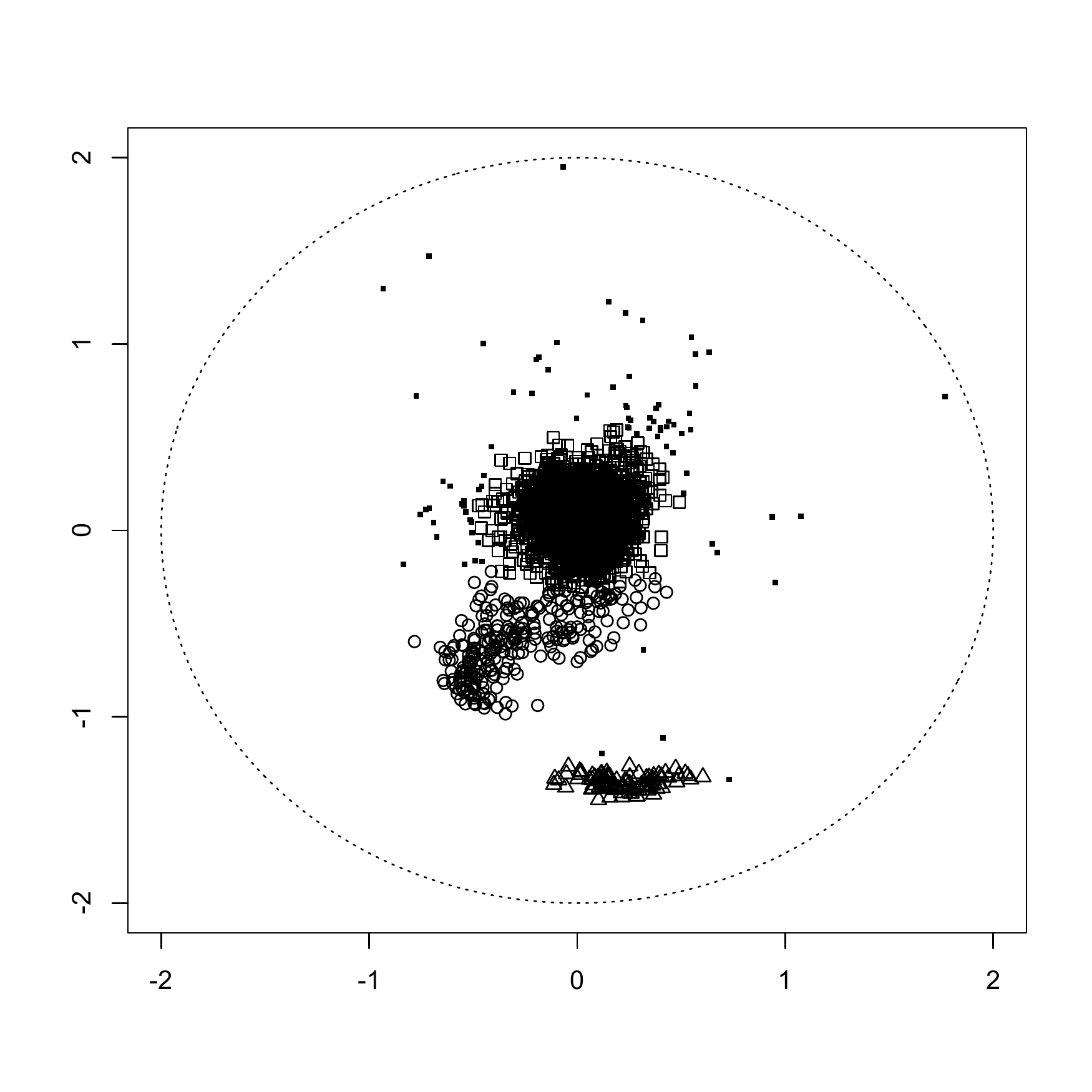}
 \caption{The helix-helix data plotted according to the
   classifications given by the `optimum' model.  Points are indicated
   by $\bigcirc$, $\square$, and $\triangle$ for the components $d$,
   $h$, and $i$ respectively.  Points for the uniform component are
   indicated by small dots.}
\label{fig:hh_fitF}
\end{figure}

Choosing the number $g$ of Kent components in our mixture model is
complicated by the fact that the fitted model found by the EM
algorithm depends heavily on the initial selections for the mean
directions.   We conduct model
selection in a stepwise manner.  In the first step a mixture model of
one Kent and one uniform component is fitted to the data using each
starting point in turn and the one with the lowest AIC chosen as the
preferred one-component model.  At each subsequent step,
$g=2, \ldots, 9$, the model with the lowest AIC in the previous step
is used as the basis for the next model and then each remaining
starting point added in turn to give a candidate model with $g$ Kent
components and one uniform component.  The $g$-component model with
the lowest AIC is retained. This process is repeated until the
nine-component model using all starting points is fitted to the data.
The final `optimal' model is chosen by comparing the `best' models
from each step and choosing the one with the lowest AIC.  Overall the
AIC is minimized by a model consisting of the three Kent components
labelled $h,d,$ and $i$ in Figure \ref{fig:hh}(a).  Component $h$ was
fitted first, $d$ second and $i$ third.  The parameter estimates for
this model are shown in the Table \ref{tab:hh_fitF}.  Figure
\ref{fig:hh_fitF} shows the allocation of data points to model
components achieved when this model is fitted to the data.  The data
are plotted with circles, squares, and triangles for components $d$,
$h$, and $i$ respectively, while observations allocated to the uniform
component are plotted with small points.

\begin{table}[tbh]
  \caption{Parameter estimates (top) and $\Delta_L$ distances (bottom)
    for the three components fitted to the helix-helix
    data ($h$, $d$, and $i$).} \label{tab:hh_fitF}
\begin{center}
\begin{tabular}{cr@{.}lr@{.}lr@{.}lr@{.}l}
\toprule
 & \multicolumn{2}{c}{uniform} & \multicolumn{2}{c}{$h$}&
   \multicolumn{2}{c}{$d$} & \multicolumn{2}{c}{$i$} \\ 
\midrule
$\kappahat$& \multicolumn{2}{c}{-}&  63&93 & 35&46& 264&34\\
$\betahat$&\multicolumn{2}{c}{-}& 1&74& 13&89 & 98&07\\
$\pihat_i$ &0&03& 0&87& 0&08& 0&03\\ 
$\gammahat_{\/1}$ & \multicolumn{2}{c}{-} & 
\multicolumn{2}{c}{$\left[\begin{array}{:}-0.23\\0.97\\0.01 \end{array}\right]$}& 
\multicolumn{2}{c}{$\left[\begin{array}{:}0.33\\-0.93\\0.14 \end{array}\right]$}& 
\multicolumn{2}{c}{$\left[\begin{array}{:}0.93\\-0.37\\0.06 \end{array}\right]$}\\
\bottomrule 
\end{tabular}

\vspace*{5mm}
\begin{tabular}{ccrrr@{\hspace{4mm}}r}
\toprule
$\Delta_L$ & \multicolumn{1}{c}{uniform} & \multicolumn{1}{c}{$h$}&
  \multicolumn{1}{c}{$d$} & \multicolumn{1}{c}{$i$} & Total \\
\midrule
2 & ~8 &    0 &   0 & 94 & 102\\
3 & 11 &   10 & 244 &  0 & 263\\
4 & 43 & 3170 &  22 &  0 & 3235\\
5+ & 19 &    6 &   0 &  0 & 25\\
\midrule
Total & 81 & 3186 & 266 & 94 & 3627\\
\bottomrule
\end{tabular}
\end{center}
\end{table}

In this case, the mixture model has recovered the physically most
plausible hydrogen bonds (those with $\Delta_L=4$) as mainly being in
component $h$, despite $\Delta_L$ not being used in the modelling
process.  Components $d$ and $i$ are mainly observations with
$\Delta_L$ being 3 and 2, respectively.  This suggests that the
mixture modelling has effectively separated the true hydrogen bonds
from atypical and/or spurious ones in the data set. The adjusted Rand
index \citep{Hubert1985} comparing the true $\Delta_L$ and allocated
component membership is 0.879, indicating a close agreement between the
two classifications.  Note that the condition numbers given
by~(\ref{eq:eccentricity}) for $d$ and $i$ are very high at 8.1 and
6.7 respectively.

%%%%%%%%%%%%%%%%%%%%%%%%%%%%%%%%%%%%%%%%%%%%%

\paragraph{Acknowledgments:} PMB was funded by a DTG award from the UK
Engineering and Physical Sciences Research Council.  The authors wish
to thank G.J.\ McLachlan for helpful discussions, and KVM thanks the
Leverhulme Trust for an Emeritus Fellowship grant.

%%%%%%%%%%%%%%%%%%%%%%%%%%%%%%%%%%%%%%%%%%%%%

\bibliographystyle{dcu}
\bibliography{refs}
\end{document}